\documentclass[11pt, reqno]{amsart}

\usepackage{etoolbox}

\newbool{bForSubmission}
\booltrue{bForSubmission}

\newbool{bIncludeStatusUpdates}

\newbool{bIncludeTOC}
\booltrue{bIncludeTOC}

\newbool{bExperimental}

\newbool{bForUs}

\usepackage{amscd,amsmath}
\usepackage{mathrsfs}
\usepackage{amsfonts}

\usepackage{amssymb, bm, xspace}
\usepackage{enumerate}

\usepackage{cite}

\usepackage[T2A,T1]{fontenc}
\usepackage[utf8]{inputenc}  
\usepackage[russian,english]{babel}
\DeclareFontEncoding{OT2}{}{}
\DeclareFontSubstitution{OT2}{cmr}{m}{n}

\usepackage{datetime}

\usepackage{enumitem}

\usepackage[dvipsnames]{xcolor}

\usepackage[pagebackref=true, colorlinks=true, citecolor=blue]{hyperref}

\usepackage[margin=1.2in]{geometry}

\usepackage{mathtools}

\mathtoolsset{centercolon}

\usepackage{cleveref}

\usepackage{stackengine}
\stackMath

\usepackage{tcolorbox}
\usepackage{tikz}
\tcbuselibrary{skins}
\usetikzlibrary{shadings}
\tcbset{
    myimage/.style={
        enhanced,
        overlay={
            \begin{scope}[shift={([xshift=1mm, yshift=7mm]frame.north west)}]
            \end{scope}}}}

\tcbset{
    skin=enhanced,
    fonttitle=\bfseries,
    interior style={white},
    segmentation style={black,solid,opacity=0.2,line width=1pt}}

\newtcolorbox{TitledBox}[2][]{
    myimage,              
    coltitle=black,       
    colbacktitle=white,   
    title=My title,
    attach boxed title to top center={
        yshift=-3mm,
        yshifttext=-1mm},
    attach boxed title to top left={
        xshift=1cm,
        yshift=-2mm},
    boxed title style={
        size=small},
    title={#2},#1}

\newcommand{\key}[1]{{\color{red}#1}}

\newcommand{\Comment}[1]{{\color{Brown}#1}}
\newcommand{\OptionalDetails}[1]{
    \ifbool{bForSubmission}
        {%
        }%
        {\begin{quote}\Comment{\footnotesize
        \medskip
        
        \noindent\textbf{Details not for submission}: \\
        \noindent#1}
        \end{quote}
        }
    }

\newbool{HaveBBM}
\booltrue{HaveBBM}

\ifbool{HaveBBM}{
	\usepackage{bbm}
    }
    {
    }

\newbool{arXivFormat}
\boolfalse{arXivFormat}
    
\newcommand{\IfarXivElse}[2]{
    \ifbool{arXivFormat}
        {#1}{#2}
    }

\renewcommand{\mathbf}[1]{\bm{#1} \textbf{ *** Use bm instead of mathbf ***}}

\newcommand{\eqn}{\begin{eqnarray}}
\newcommand{\een}{\end{eqnarray}}

\newtheorem{theorem}{Theorem}[section]
\newtheorem*{theorem*}{Theorem}				
\newtheorem{prop}[theorem]{Proposition}
\newtheorem{lemma}[theorem]{Lemma}

\newtheorem{remark}[theorem]{Remark}
\newtheorem*{remark*}{Remark}

\numberwithin{equation}{section}

\newcommand{\BoldTau}
    {{\mbox{\boldmath $\tau$}}}

\newcommand{\BB}[1]{\ensuremath{\mathbb{#1}}}
\newcommand{\R}{{\ensuremath{\BB{R}}}}

\newcommand{\iny}{{\infty}}
\newcommand{\grad}{\ensuremath{\nabla}}

\newcommand{\CharFunc}{
    \ifbool{HaveBBM}{
        \ensuremath{\mathbbm{1}}
        }
        {
        \ensuremath{\bm{1}}
        }
    }

\DeclareMathOperator{\dv}{div} %
\DeclareMathOperator{\curl}{curl} %

\DeclareMathOperator{\cond}{cond} %

\newcommand{\prt}{\ensuremath{\partial}}
\newcommand{\brac}[1]{\ensuremath{\left[ #1 \right]}}

\newcommand{\pr}[1]{\ensuremath{\left( #1 \right)}}

\DeclarePairedDelimiter{\set}{\{}{\}}

\DeclarePairedDelimiterX{\norm}[1]{\lVert}{\rVert}{#1}

\DeclarePairedDelimiterX{\abs}[1]{\lvert}{\rvert}{#1}

\newcommand\tenq[2][1]{%
	\def\useanchorwidth{T}%
	\ifnum#1>1%
		\stackunder[0pt]{\tenq[\numexpr#1-1\relax]{#2}}{\scriptscriptstyle\sim}%
	\else%
		\stackunder[1pt]{#2}{\scriptscriptstyle\sim}%
	\fi%
	}

\newcommand{\n}{{\bm{n}}}

\renewcommand{\epsilon}{\varepsilon}

\newcommand{\al}{\ensuremath{\alpha}}

\newcommand{\pdx}[2]{\frac{\prt #1}{\prt #2}}

\newcommand{\pdxn}[3]{\frac{\prt^{#3} #1}{\prt #2^{#3}}}

\newcommand{\ol}{\overline}

\renewcommand{\leq}{\leqslant}

\renewcommand{\le}{\leqslant}
\renewcommand{\ge}{\geqslant}

\newcommand{\Holders}
    {H\"{o}lder's\xspace}
    
\newcommand{\Gronwalls}
    {Gr\"{o}nwall's\xspace}

\newcommand{\Ignore}[1]{}

\newcommand{\Note}[2]{
\begin{TitledBox}{\key{#1}}
	#2
\end{TitledBox}
}

\newcommand{\PickUpHere}[1]{\Note{***}{Pick up here}}

\newcommand{\Experimental}[1]%
{%
\ifbool{bExperimental}%
	{\bigskip
	\noindent
	\textbf{\color{Brown}*** Start experimental}\\
	#1
}
{
}
}

\definecolor{Correction}{named}{red}

\crefname{cor}{Corollary}{Corollaries} 
									   
\crefname{lemma}{Lemma}{Lemmas}	       

\crefname{section}{Section}{Sections}
\Crefname{section}{Section}{Sections}

\crefname{appendix}{Appendix}{Appendices}
\Crefname{appendix}{Appendix}{Appendices}

\crefname{theorem}{Theorem}{Theorems}
\Crefname{theorem}{Theorem}{Theorems}

\crefname{prop}{Proposition}{Propositions}
\Crefname{prop}{Proposition}{Propositions}

\crefname{conj}{Conjecture}{Conjectures}
\Crefname{conj}{Conjecture}{Conjectures}

\crefname{definition}{Definition}{Definitions}
\Crefname{definition}{Definition}{Definitions}

\crefname{remark}{Remark}{Remarks}
\Crefname{remark}{Remark}{Remarks}

\crefname{assumption}{Assumption}{Assumptions}
\Crefname{assumption}{Assumption}{Assumptions}

\creflabelformat{enumi}{(#1)}
\creflabelformat{enumii}{(#1)}

\crefformat{equation}{(#2#1#3)}
\crefrangeformat{equation}{(#3#1#4) through (#5#2#6)}
\crefmultiformat{equation}
    {(#2#1#3)}%
    { and~(#2#1#3)}
    {, (#2#1#3)}
    { and~(#2#1#3)}

\newcommand{\e}{\bm{\mathrm{e}}}
\renewcommand{\f}{\bm{\mathrm{f}}}  
\renewcommand{\f}{{f}}
\newcommand{\of}{{\ol{f}}}

\newcommand{\uu}{u}

\newcommand{\vv}{v}

\newcommand{\x}{{\bm{\mathrm{x}}}}

\newcommand{\w}{w}
\newcommand{\z}{z}

\newcommand{\ou}{\ol{u}}
\newcommand{\ov}{\ol{v}}
\newcommand{\ow}{\widetilde{w}}
\newcommand{\oz}{\widetilde{z}}

\newcommand{\ovc}{{\ol{v}}}
\newcommand{\ozc}{{\widetilde{z}}}

\newcommand{\UU}{{\bm{\mathrm{U}}}}

\newcommand{\uSolSpaceN}{{C_\sigma^{N + 1, \al}(Q)}}

\newcommand{\uInitSpaceN}{{C_\sigma^{N + 1, \al}(\Omega)}}

\newcommand{\USpaceNInf}{{C_\sigma^{N + 2, \al}((0, \iny) \times \Omega)}}

\newcommand{\fspaceNInfChannel}{{C^{N + 1, \al}((0, \iny) \times \Omega) \cap C([0, \iny); H)}}

\newcommand{\wprt}{\widetilde{\prt}_t}

\newcommand{\PullMarginsIn}[1]%
{%
\begin{minipage}[c]{0.27\textwidth}%
\phantom{x}%
\end{minipage}%
\begin{minipage}[c]{0.46\textwidth}%
\begin{center}
#1%
\end{center}
\end{minipage}%
\begin{minipage}[c]{0.27\textwidth}%
\phantom{x}%
\end{minipage}%
}

\allowdisplaybreaks



\begin{document}
\newdateformat{mydate}{\THEDAY~\monthname~\THEYEAR}

\title
	[On vanishing viscosity with inflow, outflow]
	{On vanishing viscosity with inflow, outflow}

\author{Michael A. Gulas}
\address{Department of Mathematics, University of California, Riverside}
\curraddr{}
\email{mgula002@ucr.edu}

\author{James P. Kelliher}
\address{Department of Mathematics, University of California, Riverside}
\curraddr{}
\email{kelliher@math.ucr.edu}

\begin{abstract}
	We establish convergence as the viscosity vanishes 
    of solutions
	of the Navier-Stokes equations to a solution of the Euler
	equations for inflow, outflow boundary conditions. We
	extend the approach of Temam and Wang 2002, allowing
	the tangential component on outflow to be nonzero.
\end{abstract}

\maketitle

\tableofcontents


\vspace{-0.8em}

\begin{center}
Compiled on {\dayofweekname{\day}{\month}{\year} \mydate\today} at \currenttime\xspace Pacific Time
		
\end{center}

\bigskip

\noindent

\begin{quote}
\textit{We humbly dedicate this work to Roger Temam, from whom so many of us have learned so much, in celebration of his $85^{th}$ birthday.}
\end{quote}

\section{Background}

\noindent In \cite{TW02-1}, Roger Temam and Xiaoming Wang obtain a boundary layer expansion for the solution to the Navier-Stokes equations in a 3D channel periodic in the horizontal directions, with constant vertical inflow on the upper and constant vertical outflow on the lower boundary. In \cite{GHT11}, this result was extended to a three-dimensional curved domain with more general values for inflow, outflow.

In both of these references, the tangential component of the outflow velocity was required to vanish. In this paper, we remove this restriction.

First, let us define the problem we address more fully. We consider the Euler equations $(E)$ and the Navier-Stokes equations $(NS)$ with \textit{inflow, outflow} boundary conditions in a  domain $\Omega$ in $\R^d$ for $d = 2$ or $3$, either a bounded domain with boundary $\Gamma = \prt \Omega$ or a periodic channel, to which we will ultimately restrict our attention. Inflow will occur on a portion of the boundary, $\Gamma_+$, and outflow on another portion of the boundary, $\Gamma_-$ (we discuss restrictions on the portions below).

Letting $Q := (0, T) \times \Omega$, we write these systems of equations as
\begin{align*}
	(NS) \quad
	&\begin{cases}
		\prt_t \uu + \uu \cdot \grad \uu + \grad p
			= \nu \Delta \uu + \f
			& \text{in } Q, \\
		\dv \uu = 0
			& \text{in } Q, \\
		\uu(0) = \uu_0
			&\text{on } \Omega, \\
		\uu = \UU
			& \text{on } [0, T] \times \Gamma,
	\end{cases} \\
	(E) \quad
	&\begin{cases}
		\prt_t \ou + \ou \cdot \grad \ou + \grad \ol{p} = \of
            \phantom{+ \nu \Delta \uu} \,\;
			& \text{in } Q, \\
		\dv \ou = 0
			& \text{in } Q, \\
		\ou(0) = \ou_0
			&\text{on } \Omega, \\
		\ou \cdot \n = U^\n
			& \text{on } [0, T] \times \Gamma, \\
		\ou = \UU
			& \text{on } [0, T] \times \Gamma_+.
	\end{cases}
\end{align*}
The initial velocities $\uu_0$, $\ou_0$, and external forces $\f$, $\of$, are given and we will assume to be the same for both systems of equations. The vector field $\UU$ is also given, and enters into the equations only through its values on the boundary. The vector field $\n$ is the outward unit normal vector and for any vector field $\vv$ on $\Gamma$, $v^n := \vv \cdot \n$ and $\vv^\BoldTau := [\vv]_{tan}$, the tangential component of $\vv$. We will also write $\vv^\n$ for the outward normal component of $\vv$ and, when working in 2D, we define the unit tangential vector field $\BoldTau$ so that $(\n, \BoldTau)$ is in the standard orientation of $(\e_1, \e_2)$.

For $(NS)$, the full velocity on the boundary is set equal to $\UU$ on the boundary. For $(E)$, only the normal component of the velocity is set equal to $U^\n = \UU \cdot \n$ on the full boundary, though on the inflow boundary, $\Gamma_+$, the full velocity is prescribed.

There are requirements on $\UU$. To explain these, we partition the (sufficiently smooth) boundary $\Gamma = \prt \Omega$ into three portions, $\Gamma_+$, $\Gamma_-$, and $\Gamma_0$, corresponding to inflow, outflow, and impermeability, respectively. Each portion consists of a finite number of components. We fix the vector field $\UU$ on $[0, T] \times \Gamma$ and assume that
\begin{align}\label{e:UConds}
	U^\n < 0 \text{ on } \Gamma_+, \qquad
	U^\n > 0 \text{ on } \Gamma_-, \qquad
	U^\n = 0 \text{ on } \Gamma_0.
\end{align}

We require $\dv \UU = 0$, which imposes on $\UU$ the constraint that
$ 
	\int_{\Gamma_+} U^\n = - \int_{\Gamma_-} U^\n
$.

To obtain well-posedness of solutions to $(E)$, $\Gamma_0$ can be nonempty (in the classical case it is the full boundary), but to obtain the vanishing viscosity limit ($(VV)$, defined below), we require it to be empty.

To streamline the presentation and highlight the essential elements, we make three simplifying assumptions:
\begin{enumerate}
	\item
		We work in a 2D rather than 3D channel, setting
		\begin{align}\label{e:Omega}
			\Omega := [0, L] \times (0, h), \text{ periodic in } x_1.
		\end{align}
		
	\item
		We assume that the inflow and outflow velocities are
		the same and time-independent, which allows us to use a
		constant background flow,
		\begin{align}\label{e:UU}
			\UU = (a, -U) \text{ for some } a \in \R, U > 0.
		\end{align}
		
	\item
		We restrict ourselves to only showing that the vanishing
		viscosity limit holds, leaving a more complete boundary
		layer expansion to future work.
\end{enumerate}

The essential points of our proofs, we note, apply in the settings of both \cite{TW02-1} and \cite{GHT11}, as well as to time-varying inflow, outflow velocities.

Let
\begin{align*}
	H &:= \set{\uu \in (L^2(\Omega))^d
		\colon \dv \uu = 0, \uu \cdot \n = 0}, \\
	V &:= \set{\uu \in (H^1_0(\Omega))^d
		\colon \dv \uu = 0}.
\end{align*}
The global-in-time well-posedness of $(NS)$ follows from a perturbation of the classical well-posedness for no-slip boundary conditions, leading to the following (cf., Theorem III.3.2 of \cite{T2001}):

\begin{prop}\label{P:NS}
	Let $\uu_0 \in H + \UU$. In 2D, there exists a unique weak solution to $(NS)$
	in $C(0, T; H + \UU) \cap L^2(0, T; V + \UU)$.
\end{prop}

The well-posedness for the Euler equations, however, is a far more complicated issue. It requires compatibility conditions on the initial data to obtain unique short-time classical solutions with $\uu(t)$ having $C^{k, \al}(\Omega)$-regularity, $k \ge 1$. Such a need was alluded to in \cite{TW02-1} and \cite{GHT11}, and were obtained for $k = 1$ in \cite{AKMRussianOriginal} (English translation in \cite{AKM}). As in \cite{TW02-1,GHT11}, however, we will require $k > 1$.
The required compatibility conditions for higher regularity solutions were only recently obtained in \cite{GKM2023,GKM2025}.

The basic question we explore is whether or not the \textit{(classical) vanishing viscosity} limit $(VV)$ holds; that is, whether
\begin{align*}
	(VV)
		&\qquad
		\uu \to \ou \text{ in } L^\iny([0, T]; L^2(\Omega))
			\text{ as } \nu \to 0,
\end{align*}
meaning
\begin{align*}
	\sup_{t \in [0, T]} \int_\Omega \abs{\uu(t, x) - \ou(t, x)}^2 \, dx \to 0
		\text{ as } \nu \to 0.
\end{align*}

When $\UU \equiv 0$ (so $\Gamma_0 = \Gamma$), we obtain impermeable boundary conditions for the Euler equations and no-slip boundary conditions for the Navier-Stokes equations. This is the classical situation, and whether or not $(VV)$ holds in the general case has been a wide open problem for many decades. It was first shown in \cite{TW02-1}, however, that $(VV)$ holds for nontrivial $\UU$ (with $\Gamma_0 = \emptyset$), and was shown in \cite{GHT11} with $\UU$ satisfying no further conditions than we have imposed---with the critical exception, as also applies to \cite{TW02-1}, that $[\UU]_{tan} = 0$ on $\Gamma_-$. (In fact, \cite{TW02-1,GHT11} show a fair bit more than just that $(VV)$ holds, giving a more in-depth characterization of the behavior of $\uu$ near the boundary as the viscosity vanishes.)

\medskip
\noindent\textbf{Well-posedness of Euler.}
In \cite{GKM2023,GKM2025}, existence for finite time of solutions to $(E)$ along with uniqueness is proven for a multiply connected bounded domain in $3D$. To state the assumptions,

fix an integer $N \ge 0$. We say that the data has regularity $N$ for an integer $N \ge 0$ if
\begin{itemize}	
    \item
        $\Gamma = \prt \Omega$ is $C^{N + 2, \al}$ regular; 
	\item
		$\of \in \fspaceNInfChannel$;		
	\item
		$\UU \in \USpaceNInf$, $\dv \UU = 0$, and \cref{e:UConds} holds;
		
	\item
		$U_{min} := \min \set{\abs{U^\n(t, \x)} \colon (t, \x)
			\in [0, \iny) \times \Gamma_+} > 0$;

	\item
		$\ou_0 \in \uInitSpaceN$, $\ou_0^\BoldTau = \UU_0^\BoldTau$ on $\Gamma_+$,
\end{itemize}
where
\begin{align}\label{e:PrimarySpaces}
	\begin{split}
	\uInitSpaceN
		&:= \set{\ou \in C^{N + 1, \al}(\Omega) \colon \dv \ou = 0, \ou \cdot \n = U^\n(0)
					\text{ on } \Gamma}, \\
	\uSolSpaceN
		&:= \set{\ou \in C^{N, \al}(Q) \colon \curl \ou \in C^{N, \al}(Q), \,
                        \prt_t^{N + 1} \ou \in L^\iny([0, T]; C^\al(\Omega)), \\
                &\qquad
                        \dv \ou = 0, \ou \cdot \n = U^\n
					\text{ on } [0, T] \times \Gamma}, \\
        \norm{\ou}_\uSolSpaceN
            &:= \norm{\ou}_{C^{N, \al}(Q)}
                + \norm{\curl \ou}_{C^{N, \al}(Q)}
                + \norm{\prt_t^{N + 1} \ou}_{L^\iny([0, T]; C^\al(\Omega))}.
	\end{split}
\end{align}

We have the following result from \cite{GKM2025}:

\begin{prop}\label{P:GKM}[\cite{AKM,GKM2025}]
	Let $\ou_0 \in \uInitSpaceN$ and assume that the data
	has regularity $N$ and
	the compatibility conditions described in
	\cref{A:CompConds} are satisfied.
	There is a $T > 0$ such that there exists a solution
	$(\ou, \ol{p})$ to $(E)$ with
	$\ou \in \uSolSpaceN$ and $\grad \ol{p}$ in $L^\iny([0, T];
	C^{N, \al}(\Omega))$,
	which is unique up to an additive constant for the pressure.
	If $N \ge 1$, $\grad \ol{p}$ is also in $C^{N - 1, \al}(Q)$.	
\end{prop}

\begin{remark}\label{R:CompatConds}
	In particular, the compatibility conditions described in
	\cref{A:CompConds} along with the definition of the space
    $\uInitSpaceN$ require that
	$\ou_0 = \UU(0)$ on $\Gamma_+$ and
    $\ou_0\cdot \n = \UU(0) \cdot \n$ on all of $\Gamma$.
    We will impose the stronger condition that
    $\ou_0 = \UU(0)$ on $\Gamma_-$ as well.
\end{remark}

It is shown in \cite{KOS} that rather than assuming compatibility conditions, one can instead assume that the initial data is analytic, and obtain well-posedness.

Our main result is the following:

\begin{theorem}\label{T:MainResult}
    Let $\Omega$ be a 2D channel as in \cref{e:Omega} and
	$\UU$ be as in \cref{e:UU}.
	Make the assumptions as in \cref{P:GKM} with $N \ge 1$
    and, further, assume that $f = \of$ on $Q$,
    $u_0 = \ou_0$ on $\Omega$, and
    $u_0|_\Gamma = \UU(0)|_\Gamma$.
    Let $\ou$ be the
	solution to $(E)$ given by \cref{P:GKM}. Let $\uu = \uu^\nu$ be
	the solution to $(NS)$ given by \cref{P:NS}.
	There exists $T_0 \in (0, T)$ such that for all
	$t \in [0, T_0]$ and $\nu \le 1$,
	\begin{align}\label{e:MainLimit}
		\norm{\uu(t) - \ou(t)}_{L^2(\Omega)}
			&\le C (\nu t)^\frac{1}{2}e^{\frac{Ct}{2}}
			+ C \nu^{\frac{1}{2}} t.
	\end{align}
	The constant $C$ depends upon $\uu_0$, $\UU|_\Gamma$, and $\f$.
\end{theorem}

\cref{T:MainResult} improves (for a 2D channel) the results in \cite{TW02-1,GHT11} by allowing $a = [\UU]_{tan} \not\equiv 0$ on the outflow boundary, $\Gamma_-$.

The remainder of this paper is organized as follows: In \cref{S:TWProofVanishingTangentail} we give the proof of \cref{T:MainResult} in the special case in which the tangential component on outflow vanishes (so $a = 0$). This argument is that of \cite{TW02-1} when specialized to 2D and to only obtaining the vanishing viscosity limit. We extend the result to allow a nonzero tangential component on outflow in \cref{S:TangentialComponent}---this is the main novelty of this paper. In the appendix, we describe the compatibility conditions from \cite{GKM2025} required to obtain well-posedness of solutions to the Euler equations with inflow, outflow boundary conditions in the special case of \cref{e:Omega}.

\section{Zero tangential component on outflow}\label{S:TWProofVanishingTangentail}

\noindent In this section, we assume that the tangential component of the outflow velocity vanishes. Since, for simplicity, we are assuming that $\UU$ is constant on $\Omega$, this means that
\begin{align*}
	a = 0.
\end{align*}
Thus, there is also no tangential component on inflow, but we note that the value on inflow plays no significant role in the analysis.

To start, we ``homogenize'' the equations by setting
\begin{align*}
	\ov := \ou - \UU \text{ and }
	v := u - \UU
\end{align*}
in $(NS)$ and $(E)$, yielding (using that $u_0 = \ou_0$ and $\f = \of$) 
\begin{align*}
	(NS_h) \quad
     &\begin{cases} 
     \prt_t v+ v \cdot\nabla  v+ v\cdot\nabla \UU+\UU\cdot \nabla  v+\nabla  p \\
     \hspace{0.75in}=\nu\Delta v+\nu\Delta \UU+\f-\prt_t\UU-\UU\cdot\nabla \UU\, & \text{in Q},\\
      \dv v = 0\, & \text{in Q}, \\
       v(0)=u_0-\UU(0)\, & \text{on}\,\, \Omega,\\
       v = 0\, &  \text{on}\,\, [0,T]\times \Gamma,
    \end{cases} \\
 (E_h) \quad
 &\begin{cases} 
 \prt_t\ov+\ov \cdot\nabla \ov+\ov\cdot\nabla \UU+\UU\cdot \nabla \ov+\nabla \bar p\\
 \hspace{0.75in} =\of-\prt_t \UU-\UU\cdot\nabla \UU\, & \text{in Q},\\
  \dv \ov = 0\, & \text{in Q}, \\
  \ov(0)=u_0-\UU(0)\, & \text{on}\,\, \Omega,\\
  \ov\cdot \n=0\, & \text{on}\,\, [0,T]\times \Gamma,\\
  \ov = 0\, &  \text{on}\,\, [0,T]\times \Gamma_+.
\end{cases} 
\end{align*}

Letting
\begin{align}\label{wdefn}
    \ow = v - \ov,
\end{align}
the difference of $(NS_h)$ and $(E_h)$, gives
\begin{equation}\label{e:owBasicId}
    \prt_t \ow+ \UU \cdot\nabla  \ow+ \ow\cdot\nabla \UU+\nabla  (p-\bar p)=\nu\Delta \ow+\nu\Delta \ov+\nu\Delta \UU-(v\cdot\nabla v -\ov \cdot \nabla \ov).
\end{equation}

\subsection{The corrector}

Although $\ow \cdot \n = 0$ on $\Gamma$, $\ow \ne 0$ on the outflow boundary. This makes an energy argument based directly on \cref{e:owBasicId} infeasible, because difficult-to-control boundary terms would appear when integrating by parts. To get around this difficulty, we introduce a corrector $z$, a vector field on $[0, T] \times \Omega$ which must, at a minimum, equal $-\ov$ on the outflow boundary, be divergence-free, compactly supported in a fixed layer near the outflow boundary, and satisfy $z = - \ov$ on $\Gamma$. We use $z$ to form the \textit{corrected} difference,
\begin{align}
    w=\ow - z,
\end{align}
which therefore vanishes on $\Gamma$.

We use essentially the same corrector as in \cite{TW02-1}, specialized to 2D. We let 
\begin{align}\label{e:psiAndolz}
	\begin{split}
    \psi(t, x_1, x_2)
    	&:=\ovc^1(t, x_1,0)\frac{\nu}{U}(1-e^\frac{-Ux_2}{\nu}), \\
	\oz
		&:= \grad^\perp \psi
		= \Big(-\ovc^1(x_1,0)e^\frac{-Ux_2}{\nu}, \,
				\frac{\nu}{U}(1-e^\frac{-Ux_2}{\nu})\partial_1\ovc^1(x_1,0) \Big),
	\end{split}
\end{align}
noting that $\dv \oz = 0$, while on the outflow boundary, $\psi = 0$ so
\begin{align*}
	\oz
		&= (-\prt_2 \psi, \prt_1 \psi)
		= (-\ovc^1(x_1,0), 0)
		= - \ov.
\end{align*}
Thus, $\oz$ satisfies the minimal required properties of $z$, except that it is not supported near the outflow boundary; indeed, it does not vanish even on the inflow boundary. To rectify this, we cutoff $\psi$ with a $C^\iny$ function $\phi$ that depends only upon $x_2$, with $\phi \equiv 1$ on $[0, h/4]$, $\phi \equiv 0$ on $[h/2, h]$. We then let
\begin{align*}
    \z := \nabla^\perp (\phi \psi),
\end{align*}
which now satisfies all the minimal properties of a corrector.

It remains, however, to obtain various estimates on our corrector.

Since
\begin{align}\label{e:zExpanded}
	\z
		&= \phi\grad^\perp \psi + \psi \grad^\perp \phi
		= \phi \oz + \psi \grad^\perp \phi
		= \phi \oz - \psi (\phi', \, 0)
		= (\phi \ozc^1 - \phi' \, \psi, \, \phi \ozc^2),
\end{align}
for any $p \in [1, \iny]$, using also that $\abs{\psi} \le C \nu$,
\begin{align*}
	\norm{z^1}_p
		\le \norm{\ozc^1}_p + C \nu
		\le C \norm{\ozc^1}_p, \quad
	\norm{z^2}_p
		\le \norm{\ozc^2}_p.
\end{align*}

Moreover,
\begin{align*}
	\prt_1 z^1
		&= \prt_1 (\phi \ozc^1 - \phi' \, \psi)
		= \phi \prt_1 \ozc^1 - \phi' \prt_1 \psi
		= \phi \prt_1 \ozc^1 - \phi' \ozc^2, \\
	\prt_2 z^1
		&= \prt_2 (\phi \ozc^1 - \phi' \, \psi)
		= \phi \prt_2 \ozc^1 + \phi' \ozc^1
			- \phi' \prt_2 \psi - \phi'' \psi \\
		&= \phi \prt_2 \ozc^1 + 2 \phi' \ozc^1
			- \phi'' \psi, \\
	\prt_1 z^2
		&= \phi \prt_1 \ozc^2, \quad
	\prt_2 z^2
		= \phi \prt_2 \ozc^2 + \phi' \ozc^2,
\end{align*}
or,
\begin{align}\label{e:gradz}
	\begin{split}
	\grad \z
		&=
		\begin{bmatrix}
			\prt_1 z^1 & \prt_2 z^1 \\
			\prt_1 z^2 & \prt_2 z^2
		\end{bmatrix}
		= \phi \grad \oz
			+
		\begin{bmatrix}
			0 & \phi' \ozc^1 \\
			0 & \phi' \ozc^2
		\end{bmatrix}
			-
		\begin{bmatrix}
			\prt_1 (\phi' \psi) & \prt_2 (\phi' \psi) \\
			0 & 0
		\end{bmatrix} \\
		&= \phi \grad \oz
			+
		\begin{bmatrix}
			0 & \phi' \ozc^1 \\
			0 & \phi' \ozc^2
		\end{bmatrix}
			-
		\begin{bmatrix}
			\phi' \prt_1 \psi
				& \phi' \prt_2 \psi + \psi \phi'' \\
			0 & 0
		\end{bmatrix} \\
		&= \phi \grad \oz
			+
		\begin{bmatrix}
			0 & \phi' \ozc^1 \\
			0 & \phi' \ozc^2
		\end{bmatrix}
			-
		\begin{bmatrix}
			\phi' \ozc^2
				& -\phi' \ozc^1 + \psi \phi'' \\
			0 & 0
		\end{bmatrix} \\
		&= \phi \grad \oz
			+
		\begin{bmatrix}
			-\phi' \ozc^2
				& 2 \phi' \ozc^1 - \psi \phi'' \\
			0 & \phi' \ozc^2
		\end{bmatrix}.
	\end{split}
\end{align}

A direct calculation gives,
\begin{align}\label{e:gradolz}
	\nabla \oz
		&=
		\begin{bmatrix}
			-\prt_1\ov^1(x_1,0) e^\frac{-Ux_2}{\nu}
				& \ov^1 (x_1,0) \frac{U}{\nu}e^\frac{-Ux_2}{\nu} \\
			\frac{\nu}{U}(1-e^\frac{-Ux_2}{\nu})\prt_1^2\ov^1(x_1,0)
				& e^\frac{-Ux_2}{\nu}\prt_1\ov^1(x_1,0)
	\end{bmatrix}.
\end{align}

The expressions in \cref{e:psiAndolz,e:zExpanded,e:gradz,e:gradolz} lead to the following estimates (identical estimates hold for $\oz$):
\begin{align}\label{e:zBounds}
	\begin{array}{llll}
		\norm{z^1}_2\leq C\nu^\frac{1}{2},
		&\norm{z^2}_2\leq C\nu,
		&\norm{z}_2\leq C\nu^\frac{1}{2},
		&\norm{\prt_t z}_2\leq C\nu^\frac{1}{2}, \\
		\norm{\prt_1 z^1}_2\leq C\nu^\frac{1}{2},
		&\norm{\prt_2 z^1}_2\leq C\nu^{-\frac{1}{2}},
		&\norm{\prt_1 z^2}_2\leq C\nu,
		&\norm{\prt_2 z^2}_2\leq C\nu^\frac{1}{2}, \\
		\norm{z \cdot \grad z}_2\leq C\nu^\frac{1}{2}.
	\end{array}
\end{align}

Moreover, because $\ov|_{\Gamma_-} = 0$ at $t = 0$ so also $\prt_1 \ov|_{\Gamma_-} = 0$, and $\ov \in L^\iny([0, T]; C^2(\Omega))$,
\begin{align*}
	\abs{\ovc^1(t, x_1,0)}
		&\le \norm{\prt_t \ovc^1(x_1,0)}
			_{L^\iny(0, T] \times \Gamma_-)} t
		\le C t, \\
	\abs{\prt_1 \ovc^1(t, x_1,0)}
		&\le \norm{\prt_t \prt_1 \ovc^1(x_1,0)}
			_{L^\iny(0, T] \times \Gamma_-)} t
		\le C t.
\end{align*}
As a consequence, all the bounds in \cref{e:zBounds} can be improved by a factor of $t$. Only in one instance, however, will this factor improve (up to a constant factor) the rate of convergence in $(VV)$, where we will use that
\begin{align}\label{e:z2Bound}
	\norm{z}_2 \le C \nu^{\frac{1}{2}} t.
\end{align}

\begin{lemma}\label{L:prt2z1Bounds}
    There exists a constant $C$ depending only upon $U$, such that for all $\nu \le 1$,
    \begin{align*}
	   \norm{x_2^2 \prt_2 z^1}_\infty
            &\le C \nu \norm{\ovc^1(t)}_{L^\iny(\Gamma_-)}, \quad
	   \norm{x_2 \prt_2 z^1}_2
            \le C \nu^\frac{1}{2} \norm{\ovc^1(t)}_{L^\iny(\Gamma_-)}.
    \end{align*}  
\end{lemma}
\begin{proof}
From \cref{e:psiAndolz,e:gradz},
\begin{align*}
    \prt_2 z^1
        &= \phi \ov^1 (x_1,0) \frac{U}{\nu}e^\frac{-Ux_2}{\nu}
            - 2 \phi' \ovc^1(x_1,0)e^\frac{-Ux_2}{\nu}
            - \psi \phi''.
\end{align*}
Along with the simple bound, $\abs{\psi} \le C \norm{\ovc^1(t)}_{L^\iny(\Gamma_-)} \nu$, we have
\begin{align*}
	\norm{x_2^2 \prt_2 z^1}_\infty
		&\le \norm{\ovc^1(t)}_{L^\iny(\Gamma_-)} \,
            \brac{
			\norm{x_2^2 (U/\nu)
				e^{\frac{-Ux_2}{\nu}}}_\infty
            + C \norm{x_2^2 e^{\frac{-Ux_2}{\nu}}}_\iny
            + C \nu} \\
        &\le C (\nu + \nu^2) \norm{\ovc^1(t)}_{L^\iny(\Gamma_-)}.
\end{align*}
We used that

$\norm{x^2 c e^{-cx}}_\infty = 4 e^{-2}/c$, so also $\norm{x^2 e^{-cx}}_\infty = 4 e^{-2}/c^2$.

This gives the first bound. For the second bound,
\begin{align*}
    \norm{x_2 \prt_2 z^1}_2^2
        &\le 3 \norm{\ovc^1(t)}_{L^\iny(\Gamma_-)}^2
            \brac{
                \frac{U^2}{\nu^2} \int_0^\iny x_2^2 e^{-\frac{2U x_2}{\nu}}
                + C
                \int_0^\iny x_2^2 e^{-\frac{2U x_2}{\nu}}
                + C \nu^2
            } \\
        &= 3 \norm{\ovc^1(t)}_{L^\iny(\Gamma_-)}^2
            \brac{
                \pr{\frac{U^2}{\nu^2} + C}
                    \frac{\nu^3}{4 U^3}
                + C \nu^2
            } \\
        &\le C\norm{\ovc^1(t)}_{L^\iny(\Gamma_-)}^2
            (\nu + \nu^3 + \nu^2).
        \qedhere
    \end{align*}
\end{proof}

In applying the corrector estimates, we will employ the following form of Hardy's inequality  (see, for instance, Lemma II.1.10 \cite{T2001}):

\begin{lemma}\label{L:Hardy}

For any $f\in H_0^1(\Omega)$,
$
    \displaystyle \norm[\bigg]{\frac{f}{x_2}}_{L^2(\Omega)}
        \leq C_H \norm{\nabla f}_{L^2(\Omega)}
$.

\end{lemma}

\subsection{Energy argument}
Written in terms of $w$, \cref{e:owBasicId} becomes
\begin{align*}
	\prt_t  w+ &\UU \cdot\nabla  w+  w\cdot\nabla \UU+\nabla  (p-\bar p) \\
		&= \nu\Delta \ow+\nu\Delta \ov+\nu\Delta \UU-(v\cdot\nabla v -\ov \cdot \nabla \ov)-\prt_tz-\UU\cdot \nabla z-z\cdot \nabla \UU.
\end{align*}
Multiplying by $w$ and integrating over $\Omega$ yields
\begin{align*}
	\frac{1}{2} \frac{d}{dt}\norm{w}_2^2
		= &\nu(\Delta \ow, w) + \nu(\Delta \ov, w)
			+ \nu(\Delta \UU, w)-(\UU \cdot\nabla  w, w)
			-  (w\cdot\nabla \UU, w) \\
    	&- (v\cdot\nabla v -\ov \cdot \nabla \ov, w)
		- (\prt_tz,w)-(\UU\cdot \nabla z,w)
		- (z\cdot \nabla \UU,w).
\end{align*}
Integrating by parts gave $(\nabla (p-\bar p),w)=0$.

Using that $\ow=w+z$, integrating by parts gives
 
$$\nu(\Delta \ow, w)=-\nu(\nabla \ow, \nabla w)=-\nu(\nabla w, \nabla w)-\nu(\nabla z, \nabla w)$$
$$=-\nu\norm{\grad w}_2^2-\nu(\nabla z, \nabla w), $$
since $w=0$ on $\Gamma$. Bringing $-\nu\norm{\grad w}_2^2$ to the left hand side to set up an energy argument,

\begin{align}\label{EnergySetup}
	\begin{split}
	\frac{1}{2}&\frac{d}{dt}\norm{w}_2^2
			+\nu\norm{\grad w}_2^2
		= -\nu(\nabla z, \nabla w) + \nu(\Delta \ov, w)
			+ \nu(\Delta \UU, w)
			-(\UU \cdot\nabla  w, w) \\
		&-  (w\cdot\nabla \UU, w)
			- (v\cdot\nabla v -\ov \cdot \nabla \ov, w)
			- (\prt_tz,w)-(\UU\cdot \nabla z,w)
			-(z\cdot \nabla \UU,w).
	\end{split}
\end{align}
In what follows, we bound the nine terms on the right hand side of \cref{EnergySetup}.

\subsection{The easier terms}\label{S:EasierTerms}
First, integrating by parts then applying the Cauchy-Schwarz and Young's inequalities,
\begin{align}\label{E1A}
	\nu(\Delta\ov, \w) 
		&= - \nu (\grad \ov, \grad \w)
		\le \nu \norm{\grad \ov}_2 \norm{\grad \w}_2
		\le C \nu \norm{\grad \ov}_2^2
			+ \frac{\nu}{6} \norm{\grad \w}_2^2.
\end{align}
Since we must assume that $\UU \in C^2(\Omega)$, we just estimate,
\begin{align}\label{E1B}
	\nu(\Delta \UU, w)
		&\le \nu \norm{\Delta \UU}_2 \norm{\w}_2
		\le C \nu^2 + C \norm{\w}_2^2.
\end{align}

For $-(\prt_t z,w)$, we have $\norm{\prt_t z}_2\leq C\nu^\frac{1}{2}$ by \cref{e:zBounds}, so
\begin{align}\label{E2}
    -(\prt_t z, w)\leq \norm{\prt_t z}_2 \,\norm{w}_2\leq C\nu^\frac{1}{2} \norm{w}_2\leq C\nu +C\norm{w}_2^2.
\end{align}

Using that $\nabla w \cdot w = \frac{1}{2}\nabla |w|^2$, integrating by parts, we find
\begin{align}\label{E3}
	-(\UU\cdot \nabla w, w)
		&=- \frac{1}{2} \int_\Omega \UU \cdot \grad \abs{w}^2
    	=\frac{1}{2}\int_\Omega \dv \UU |w|^2
			- \frac{1}{2} \int_{\Gamma} U^n |w|^2
		= 0,
\end{align}
since $w=0$ on $\Gamma$ and $\UU$ is divergence-free.

For $(w\cdot \nabla \UU, w)$ we simply apply \Holders inequality,
\begin{align}\label{E4}
    (w\cdot \nabla \UU, w)\leq \norm{\grad \UU}_\infty \norm{w}_2^2\leq C\norm{w}_2^2.
\end{align}

The final easier term, $-(z\cdot\nabla \UU, w)=(z\cdot \nabla w, \UU)$, we split into several terms to bound individually:
$$(z\cdot \nabla w, \UU)=\int_\Omega z^1\prt_1w^1U^1 + \int_\Omega z^1\prt_1w^2U^2 + \int_\Omega z^2\prt_2w^1U^1 + \int_\Omega z^2\prt_2w^2U^2.$$
Since $U^1=0$, the first and third integral vanish, and using the product rule, we obtain
$$(z\cdot \nabla w, \UU)=\int_\Omega z^1\prt_1(w^2U^2) - \int_\Omega z^1\prt_1(U^2)w^2 + \int_\Omega z^2\prt_2(w^2U^2) - \int_\Omega z^2\prt_2(U^2)w^2. $$
From here we can integrate the first and third integrals by parts, noticing the boundary term vanishes since $w^2=0$ on $\Gamma$:
$$(z\cdot \nabla w, \UU)=-\int_\Omega \prt_1z^1(w^2U^2) - \int_\Omega z^1\prt_1(U^2)w^2 -\int_\Omega \prt_2z^2(w^2U^2) - \int_\Omega z^2\prt_2(U^2)w^2.$$
We finish, using \Holders and Young's inequalities,
\begin{align}\label{E5}
	\begin{split}
	(z\cdot \nabla w, \UU)
		&\leq \norm{\prt_1 z^1}_2\,\norm{w}_2\norm{\UU}_\infty
			+ C\nu^\frac{1}{2}\norm{w}_2
			+ \norm{\prt_2 z^2}_2\,\norm{w}_2\norm{\UU}_\infty
			+ C\nu\norm{w}_2, \\
		&\leq  C\nu^\frac{1}{2}\norm{w}_2\norm{\UU}_\infty
				+ C\nu^\frac{1}{2}\norm{w}_2
				+ C\nu^\frac{1}{2}\,\norm{w}_2\norm{\UU}_\infty
				+ C\nu\norm{w}_2 \\
    	&\leq C\nu+C\nu^2+C\norm{w}_2^2.
	\end{split}
\end{align}

\subsection{The nonlinear term}\label{S:NonlinearTerm}

For the nonlinear term, $(v\cdot\nabla v -\ov \cdot \nabla \ov, w)$, we use that $v=w+\ov + z$, so
\begin{align}\label{e:NonlinearBreakdown}
	\begin{split}
	v\cdot\nabla v - &\ov \cdot \nabla \ov
		= (w+\ov + z)\cdot\nabla (w+\ov + z)
			- \ov \cdot \nabla \ov, \\
		&= (w+\ov + z)\cdot\nabla w +(w+\ov + z) \cdot\nabla \ov
			+ (w+\ov + z\cdot)\nabla z
			- \ov \cdot \nabla \ov, \\
		&= v\cdot \nabla w + w\cdot \nabla \ov +
			\ov\cdot \nabla \ov + z\cdot \nabla \ov
			+ w\cdot \nabla z + \ov\cdot \nabla z
			+ z\cdot \nabla z -\ov \cdot \nabla \ov, \\
		&=v\cdot \nabla w + w\cdot \nabla \ov
			+ z\cdot \nabla \ov + w\cdot \nabla z
			+ \ov\cdot \nabla z + z\cdot \nabla z.
	\end{split}
\end{align}

We must bound each of the six terms on the right hand side of \cref{e:NonlinearBreakdown} when paired with $w$.

The first of the six terms from \cref{e:NonlinearBreakdown} we integrate by parts, giving
\begin{align}\label{NL1}
	\begin{split}
		(v\cdot\nabla w, w)
			&= \int_\Omega (v\cdot\nabla w)\cdot w
			= \frac{1}{2}\int_\Omega v \cdot \nabla |w|^2 \\
			&= -\frac{1}{2}\int_\Omega \dv v |w|^2
				+ \frac{1}{2} \int_{\Gamma} (v\cdot \n)|w|^2
			=0,
	\end{split}
\end{align}
since $v \cdot \n = 0$ on the boundary and $v$ is divergence-free.

For the second term from \cref{e:NonlinearBreakdown}, the (generalized) \Holders inequality gives
\begin{align}\label{NL2}
    (w\cdot\nabla \ov, w)\leq \norm{\grad \ov}_\infty\, \norm{w}_2\,\norm{w}_2\leq C\norm{w}_2^2.
\end{align}
For the third term from \cref{e:NonlinearBreakdown}, \Holders inequality, \cref{e:zBounds}, and Young's inequality give
\begin{align}\label{NL3}
    (z\cdot \nabla \ov, w)\leq \norm{\grad \ov}_\infty\, \norm{z}_2\,\norm{w}_2\leq C\nu+C\norm{w}_2^2.
\end{align}
For the sixth term from \cref{e:NonlinearBreakdown}, we use the Cauchy-Schwarz inequality, \cref{e:zBounds}, and Young's inequality, giving
\begin{align}\label{NL6}
    (z\cdot \nabla z, w)\leq \norm{z \cdot \grad z}_2\,\norm{w}_2\leq C\nu^{\frac{1}{2}}\norm{w}_2\leq C\nu + C\norm{w}_2^2.
\end{align}

The remaining two terms, the fourth and fifth, from \cref{e:NonlinearBreakdown} are more delicate. We will expand each of these terms using indices and we will find that terms involving $\prt_2z^1$ are the most delicate because, as \cref{e:zBounds} indicates, $\prt_2z^1$ introduces a factor of $\nu^{-\frac{1}{2}}$.


The fourth term we write in terms of indices and treat each term individually:
\begin{align}\label{e:FourthNonlinearTerm}
	(w\cdot \nabla z, w)
		=\int_\Omega w^i\prt_iz^jw^j
		=\int_\Omega \pr{w^1\prt_1z^1w^1 + w^1\prt_1z^2w^2
			+ w^2\prt_2z^1w^1 + w^2\prt_2z^2w^2}.
\end{align}
.
Three of these terms are easy: using \Holders and \cref{e:zBounds}, we find
$$( w^1\prt_1z^1, w^1)\leq \norm{w^1}_2 \norm{\prt_1 z^1}_\infty \norm{w^1}_2 \leq C\norm{w}_2^2,$$
$$( w^1\prt_1z^2, w^2)\leq \norm{w^1}_2 \norm{\prt_1 z^2}_\infty \norm{w^2}_2 \leq C\nu\norm{w}_2^2,$$
$$( w^2\prt_2z^2, w^2)\leq \norm{w^2}_2 \norm{\prt_2 z^2}_\infty \norm{w^2}_2 \leq C\norm{w}_2^2.$$

\begin{remark}\label{R:CauseOneofT0}
The final term in \cref{e:FourthNonlinearTerm} to bound, $\left(w^2\prt_2z^1, w^1\right)$, is problematic, and will restrict us to obtaining \cref{e:MainLimit} for a time $T_0$ possibly less than the existence time of the solution to the Euler equations.
\end{remark}

To bound $\left(w^2\prt_2z^1, w^1\right)$, we use \cref{L:Hardy} followed by \cref{L:prt2z1Bounds}, 
\begin{align*}
	\left(w^2\prt_2z^1, w^1\right)
		& =\left(\frac{w^2}{x_2} (x_2^2\prt_2z^1),
				\frac{w^1}{x_2}\right)
			\leq \norm[\bigg]{\frac{w^2}{x_2}}_2\,
			\norm{x_2^2\prt_2z^1}_\infty\,
			\norm[\bigg]{\frac{w^1}{x_2}}_2
			\\
		&\leq \norm{x_2^2 \prt_2z^1}_\infty C_H^2
			\norm{\grad w}_2^2
        \le C_H^2 C \nu \norm{\ovc^1(t)}_{L^\iny(\Gamma_-)}
            \norm{\grad w}_2^2.
\end{align*}

Now, $\ovc^1 = 0$ on $\Gamma$ at $t = 0$ (see \cref{R:CompatConds}) or $\ovc^1 (0, x_1,0) = 0$, and $\ov$ is continuous on $[0, T] \times \ol{\Omega}$. Therefore, given any $C_0 > 0$ there exists a time $T_0 > 0$ so that for all $t \in [0, T_0]$,
\begin{align}\label{timebound}
    \norm{\ovc^1(t)}_{L^\iny(\Gamma_-)} \leq \frac{1}{C_0}.
\end{align}
Hence, for $\nu \le 1$
\begin{align*}
    \left(w^2\prt_2z^1, w^1\right)\leq \frac{C_H^2}{C_0} C\nu \norm{\grad w}_2^2 \leq \frac{\nu}{6} \norm{\grad w}_2^2
\end{align*}
by choosing $C_0$ and hence $T_0$ sufficiently small. This gives 
\begin{align}\label{NL4}
    (w\cdot \nabla z, w)
    	&\leq C\norm{w}_2^2
			+C \nu \norm{w}_2^2
			+ \frac{\nu}{6} \norm{\grad w}_2^2
		\leq C\norm{w}_2^2
			+ \frac{\nu}{6} \norm{\grad w}_2^2.
\end{align}


Finally, the fifth term from \cref{e:NonlinearBreakdown} we also write in terms of indices and treat each term individually:
\begin{align*}
	(\ov \cdot \nabla z , w)
		&=\int _\Omega \ov^i \prt_i z^jw^j
		=(\ovc^1 \prt_1 z^1, w^1)
			+ (\ovc^1 \prt_1 z^2, w^2) + (\ovc^2 \prt_2 z^1, w^1)
			+ (\ovc^2 \prt_2 z^2, w^2).
\end{align*}
The first, second, and fourth terms are easily bounded:
$$(\ovc^1 \prt_1 z^1, w^1)\leq \norm{\ovc^1}_\infty\, \norm{\prt_1 z^1}_2\,\norm{w^1}_2\leq C\nu^\frac{1}{2}\,\norm{w}_2\leq C\nu +C\norm{w}_2^2,$$
$$(\ovc^1 \prt_1 z^2, w^2)\leq \norm{\ovc^1}_\infty\, \norm{\prt_1z^2}_2\,\norm{w^2}_2\leq C\nu\,\norm{w}_2 \leq C\nu^2 +C\norm{w}_2^2,$$
$$(\ovc^2 \prt_2 z^2, w^2)\leq \norm{\ovc^2}_\infty\, \norm{\prt_2 z^2}_2\,\norm{w^2}_2\leq C\nu^\frac{1}{2}\,\norm{w}_2 \leq C\nu +C\norm{w}_2^2.$$
For the third term, we use the mean value theorem in $x_2$ followed by \cref{L:prt2z1Bounds}, giving
\begin{align*}
    (\ovc^2 \prt_2 z^1, w^1)
        &= \pr{\frac{\ovc^2}{x_2} (x_2\prt_2z^1), w^1}
        \le \norm[\bigg]{\frac{\ovc^2}{x_2}}_\iny
            \norm{x_2\prt_2z^1}_2
            \norm{w}_2
        \le C C_H \norm{\prt_2 \ovc^2}_\iny \nu^\frac{1}{2}\norm{w}_2 \\
        &\le C\nu + C\norm{w}_2^2,
\end{align*}
by Young's inequality. Thus, we have
\begin{align}\label{NL5}
    (\ov \cdot \nabla z , w)\leq C\nu + C\nu^2+C\norm{w}_2^2.
\end{align}


Collecting the estimates in \cref{NL1,NL2,NL3,NL6}, \cref{NL4}, and \cref{NL5},
\begin{align}\label{NL}
	\begin{split}
	(v\cdot\nabla v - &\ov \cdot \nabla \ov, w)
		\leq 0 + C\norm{w}_2^2 + (C\nu+\norm{w}_2^2)
			+ (C\nu+C\norm{w}_2^2) \\
			&\qquad\qquad\qquad\qquad
			+(C\norm{w}_2^2
			+ \frac{\nu}{6} \norm{\grad w}_2^2)
			+ (C\nu + C\nu^2 + C\norm{w}_2^2) \\
		&\leq C\nu + C\norm{w}_2^2 + \frac{\nu}{6} \norm{\grad w}_2^2.
	\end{split}
\end{align}

\subsection{The heart of the matter} 

From \cref{EnergySetup}, using \cref{E1A,E1B,E2,E3,E4,E5,NL}, we have
\begin{align}\label{EnergySetup2}
	\begin{split}
	\frac{1}{2} &\frac{d}{dt} \norm{w}_2^2
			+\nu\norm{\grad w}_2^2 \\
		&\leq -\nu(\nabla z, \nabla w)+(C\nu^2+C\norm{w}_2^2)
			+ (C\nu+C\norm{w}_2^2) + 0 + C\norm{w}_2^2 \\
		&\qquad
		+ (C\nu+C\nu^2 + C\norm{w}_2^2)
			+ (C\nu + C\nu^2 + C\norm{w}_2^2
				+ \frac{\nu}{3}\norm{\grad w}_2^2) 
			- (\UU\cdot \nabla z,w) \\
		&\leq C\nu + C\nu^2 + C\norm{w}_2^2
				+ \frac{\nu}{3}\norm{\grad w}_2^2 
			-\nu(\nabla z, \nabla w)  - (\UU\cdot \nabla z,w).
	\end{split}
\end{align}
Two terms in \cref{EnergySetup}, $-\nu(\nabla z, \nabla w)$ and $-(\UU\cdot \nabla z,w)$, remain to be bounded. The first of these arises from the diffusive term in $(NS_h)$, while the second originates in the nonlinear terms of $(E)$ and $(NS)$, from which it was separated as a consequence of treating $(E_h)$ and $(NS_h)$ as perturbations of $(E)$ and $(NS)$.
Rather than bound these terms separately, we will bound their sum: in doing that, we will see the central reason, at the heart of the innovation in \cite{TW02-1}, that the vanishing viscosity limit can be obtained for inflow, outflow boundary conditions.

To better appreciate this innovation, let us let us briefly consider how the energy argument we are making differs from the analogous argument in the classical case in which one enforces no-slip conditions ($\uu = 0$ on $\Gamma$) on $(NS)$ and impermeable conditions ($\ou \cdot \n = 0$ on $\Gamma$) on $(E)$. (The  situation regarding the classical vanishing viscosity limit is perhaps most clearly expressed in Tosio Kato's \cite{Kato1983}.)

The classical analog of all the terms in \cref{EnergySetup} that we have so far bounded can still be bounded---or they do not appear because, in effect, $\UU = 0$. (In particular, $\left(w^2\prt_2z^1, w^1\right)$ of \cref{R:CauseOneofT0} can be controlled when
$\uu_0 = 0$ on $\Gamma$, much as was done here, for short time.) The diffusive term, $-\nu(\nabla z, \nabla w)$, however, cannot be controlled---neither in the classical nor in our setting.

The additional term, $-(\UU\cdot \nabla z,w)$, absent in the classical case, also cannot be controlled. This seeming disadvantage turns out to be an advantage, for we will find that the combined terms,
\begin{align*}
	I
		&:= - \nu(\grad \z, \grad \w) - (\UU \cdot \grad \z, \w),
\end{align*}
can, in fact, be controlled, as long as $\UU \cdot \n$ never vanishes on the outflow boundary.

The first step in bounding $I$ is to observe that $-(\UU\cdot \nabla z,w)=(\UU\cdot \nabla w,z)=(z\otimes \UU, \nabla w)$, after integrating by parts, so
\begin{equation*} 
    I=(\nu\nabla z - z\otimes \UU, \nabla w).
\end{equation*}


\subsection{Bounding the inertial plus diffusive terms}

From \cref{e:gradolz},
\begin{align}\label{IbigOnotation}
	\begin{split}
	\nu\nabla \oz - \oz\otimes \UU
		&=
		\begin{bmatrix}
			-\nu\prt_1\ovc^1(x_1,0) e^\frac{-Ux_2}{\nu} - U^1 \ozc^1
				&\nu\ovc^1 (x_1,0) \frac{U}{\nu}e^\frac{-Ux_2}{\nu}
					-U^2 \ozc^1 \\
			\frac{\nu^2}{U}(1-e^\frac{-Ux_2}{\nu})\prt_1^2\ovc^1(x_1,0)
					-U^1 \ozc^2
				&\nu e^\frac{-Ux_2}{\nu}\prt_1\ovc^1(x_1,0)
					- U^2 \ozc^2
		\end{bmatrix} \\
		&=
		\begin{bmatrix}
			O(\nu^{\frac{3}{2}}) & 0 \\
			O(\nu^2) & O(\nu)
		\end{bmatrix}
		\text{ in } L^2(\Omega),
	\end{split}
\end{align}
where we used that $U^1 = 0$ and that $U \ne 0$
\begin{align}\label{e:KeyI}
	\nu \prt_2 \ozc^1 - U^2 \ozc^1
		&= \nu\ov^1 (x_1,0) \frac{U}{\nu} e^\frac{-Ux_2}{\nu}
			- (-U)(-\ov^1(x_1,0) e^\frac{-Ux_2}{\nu})
		= 0.
\end{align}
Hence,
\begin{align*}
    \norm{\nu\nabla \oz - \oz \otimes \UU}_2
    	\leq C \nu.
\end{align*}

\begin{remark}\label{R:NeedsConstantU}
	The term $\prt_2 \ozc^1=\ov^1 (x_1,0) \frac{U}{\nu}
	e^\frac{-Ux_2}{\nu}$ cannot be controlled on its own;
	it is the cancellation by $U^2 \ozc^1$ that allows
	control. This was a virtue of the
	corrector chosen in \cite{TW02-1}.
\end{remark}

Then, from \cref{e:zExpanded,e:gradz},
\begin{align*}
	\nu \grad \z - \z \otimes \UU
		&= \phi (\nu\nabla \oz - \oz\otimes \UU)
			+
		\nu
		\begin{bmatrix}
			-\phi' \ozc^2
				& 2 \phi' \ozc^1 - \psi \phi'' \\
			0 & \phi' \ozc^2
		\end{bmatrix}
		+ \psi (\phi', \, 0) \otimes \UU.
\end{align*}

Other than the first and last terms above, the largest magnitude term derives from $\ozc^1$, which is, however, multiplied by $\nu$. Hence, all the terms beyond the first have $L^2$ norm no larger than $C \nu$, and we conclude that
\begin{align*} 
    \norm{\nu\nabla \z - \z \otimes \UU}_2
    	\leq C\nu.
\end{align*}

Therefore, applying \Holders and Young's inequality,
\begin{align}\label{I Bound with cutoff Function}
	\abs{I}
		&\le \norm{\nu\nabla \z - \z \otimes \UU}_2
			\norm{\grad \w}_2
		\le C \nu + \frac{\nu}{6} \norm{\grad \w}_2^2.
\end{align}

\subsection{Completing the energy argument}

From \cref{I Bound with cutoff Function,EnergySetup2}, we have
\begin{align*}
	\begin{split}
	\frac{1}{2}&\frac{d}{dt}\norm{w}_2^2
			+\nu\norm{\grad w}_2^2
		\leq C\nu + C\nu^2 + C\norm{w}_2^2
				+ \frac{\nu}{2}\norm{\grad w}_2^2 
	\end{split}
\end{align*}
so
\begin{align*}
	\frac{d}{dt}\norm{w}_2^2
			+ \nu\norm{\grad w}_2^2
		\leq C\nu +C\norm{w}_2^2+C\nu\norm{\grad w}_2^2,
\end{align*}
since we have assumed that $\nu \le 1 $. Then, integrating over time,

$$\norm{w(t)}_2^2+\nu \int_0^t \norm{\nabla w(s)}_2^2\,ds \leq \norm{w(0)}_2^2 +C\nu t +\int_0^t C\norm{w(s)}_2^2\,ds. $$

Because $u_0 = \ou_0$, $w(0)=0$, and \Gronwalls inequality yields
$$\norm{w(t)}_2^2 + \nu \int_0^t \norm{\nabla w(s)}_2^2 \, ds \leq C\nu te^{Ct}.$$

Thus,
$$\norm{w(t)}_2\leq C\nu^\frac{1}{2}t^\frac{1}{2}e^{\frac{Ct}{2}}.$$

By virtue of \cref{e:z2Bound}, then,
\begin{align*}
	\norm{\uu(t) - \ou(t)}_2
		&= \norm{\ow}_2
		= \norm{\w + \z}_2
		\le \norm{w}_2 + \norm{\z}_2
		\le C (\nu t)^\frac{1}{2}e^{\frac{Ct}{2}}
			+ C \nu^{\frac{1}{2}} t
\end{align*}
for all $t\in [0, T_0]$.

This completes the proof of \cref{T:MainResult} when the tangential component of the outflow velocity vanishes ($a = 0$).


\section{Nonzero tangential component on outflow}\label{S:TangentialComponent}

\noindent In this section we allow a nonzero tangential component $\UU=(a,-U)$ for some $U>0$ and $a\in \mathbb R$ in the proof of \cref{T:MainResult}.

\begin{remark}
Key to our argument is that we use exactly the same form for the corrector that we used when $a = 0$. This is because the corrector corrects for the value of $\ov$ on the boundary, which, even for $a = 0$, has a nonzero value on the outflow boundary. Hence, while $\ovc^1$ incorporates in it the value of $a$, via $\ov = \ou - \UU = \ou - (a, -U)$, the only properties of $\ov$ that we used in \cref{S:TWProofVanishingTangentail} were its regularity on $[0, T] \times \Omega$, and that does not change.
\end{remark}

We proceed by examining each of the terms in \cref{EnergySetup} and reporting on any differences. We start with the easier terms bounded in \cref{S:EasierTerms}.

There is no change in \cref{E1A}, \cref{E1B}, \cref{E3}, and \cref{E4}, but the last of these easier pieces will be slightly different now. For $-(z\cdot\nabla \UU, w)=(z\cdot \nabla w, \UU)$ we integrate by parts and find, now allowing $U^1 = a$, that
$$(z\cdot \nabla w, \UU)=\int_\Omega z^1\prt_1w^1a + \int_\Omega z^1\prt_1w^2U^2 + \int_\Omega z^2\prt_2w^1a + \int_\Omega z^2\prt_2w^2U^2.$$
There is no change to the second and fourth integrals. For the first term,
$$\int_\Omega z^1\prt_1w^1a= -a\int_\Omega \prt_1z^1w^1\leq |a| \,\norm{\prt_1 z^1}_2\,\norm{w^1}_2\leq C\nu^\frac{1}{2}\norm{w}_2\leq C\nu+C\norm{w}_2^2.$$
This followed from integration by parts and the Cauchy-Schwarz and Young's inequalities.

In a similar manner we find 
$$\int_\Omega z^2\prt_2w^1a=-a\int_\Omega \prt_2z^2w^1\leq |a|\,\norm{\prt_2 z^2}_2\,\norm{w}_2\leq C\nu^\frac{1}{2}\norm{w}_2\leq C\nu+\norm{w}_2^2.$$
Thus with nonzero tangential component we find:
\begin{align}\label{E5witha}
    -(z\cdot\nabla \UU, w)\leq C\nu+C\nu^2+C\norm{w}_2^2+(C\nu+C\norm{w}_2^2)+(C\nu+C\norm{w}_2^2).
\end{align}

There are no changes in the terms arising from the nonlinear term of \cref{S:NonlinearTerm}.

Finally for $I$, akin to notation in \cref{IbigOnotation} we find
$$\nu\nabla \oz - \oz \otimes \UU=\begin{bmatrix}
O(\nu^{\frac{3}{2}})-a\ovc^1(x_1,0) e^\frac{-Ux_2}{\nu} & 0\\
O(\nu^2)+O(\nu) & O(\nu)
\end{bmatrix}
\text{ in } L^2(\Omega).$$
We just have to handle the new exponentially decaying piece, $- a \ovc^1 e^{-\frac{Ux_2}{\nu}}$, which appears in the energy argument as
\begin{align*}
	(- a \ovc^1 e^{-\frac{Ux_2}{\nu}}, \prt_1 w^1)
        &= (a \ozc^1, \prt_1 w^1)
        = (-a \prt_1 \ozc^1, w^1).
\end{align*}
Here, we integrated by parts in $x_1$, there being no boundary term because the functions are periodic in $x_1$. Applying \Holders then Young's inequality and using \cref{e:zBounds},
\begin{align}\label{expobound}
    (a\prt_1 \ovc^1 e^{-\frac{Ux_2}{\nu}}, w^1)
        \le \abs{a} \norm{\prt_1 \ozc^1}_2 \norm{w}_2
        \leq C\nu^\frac{1}{2}\norm{w}_2\leq C\nu+C\norm{w}_2^2. 
\end{align}
We note that in obtaining this bound, the detailed structure of the corrector was not required, only the bounds in \cref{e:zBounds}.
Other than the values of constants, our bounds are the same as in \cref{S:TWProofVanishingTangentail}, so we obtain \cref{T:MainResult} for $a \ne 0$.


\section{Concluding Remarks}\label{S:Conclusion}

\noindent In this paper, we have assumed, as in \cite{TW02-1}, that the background flow $\UU$ is constant throughout $\Omega$. This derives from the assumption that the boundary conditions themselves are constant along each boundary, and hence, via the requirement in \cref{e:UConds}, $U^n|_{\Gamma_+} = - U^n|_{\Gamma_-}$ and both are constant along their respective boundary components. Because of the simple geometry, this allows for a background flow constant throughout the domain.

This requirement is dropped in \cite{GHT11}, where the boundary conditions are allowed to vary along the (curved) boundary, and leads to a background flow $\UU$ which varies over the domain. This requires an adaptation to the corrector, which nonetheless retains the same key bounds that allow the convergence in \cref{T:MainResult}.

We can gain an appreciation of why the adaptations in the \cite{GHT11} work by considering what would happen if we allowed $\UU$ to vary in our case of a 2D channel. First, observe that of all the estimates we made, it is only in $I = (\nu\nabla z - z\otimes \UU, \nabla w)$ that we used the detailed structure of the corrector $\z$, and that, only for the first row second column of $\nu\nabla z - z\otimes \UU$ and $\grad w$. For all other terms, we needed only the bounds in \cref{e:zBounds}.

So let us reexamine the key term in $I$, $(\nu \prt_2 \ozc^1 - U^2 \ozc^1, \prt_2 w^1)$. As in \cref{e:KeyI}, $\nu \prt_2 \ozc^1 - U^2 \ozc^1$ vanishes identically. This is a very special property that uses not only the special structure of our corrector, but requires that $U^2 = -U$ be constant throughout the domain. (Because $\UU$ is divergence-free, $\prt_2 U^2 = - \prt_1 U^1$, which means that both $U^1$ and $U^2$ must be constant, at least along the boundary.)

Consider what happens to this key term if $U^2$ varies over $\Omega$. Proceeding as in \cref{I Bound with cutoff Function}, this key term becomes, integrating by parts,
\begin{align*}
    (\nu \prt_2 \ozc^1 - &U^2 \ozc^1, \prt_2 w^1)
        = (\nu \prt_2^2 \ozc^1 - \prt_2 (U^2 \ozc^1), w^1) \\
        &= (\nu \prt_2^2 \ozc^1 - U^2 \prt_2 \ozc^1, w^1)
            + (\prt_2 U^2 \ozc^1, w^1)
        =: I_1 + I_2.
\end{align*}
Applying \Holders and Youngs inequalities,
\begin{align*}
    \abs{I_1}
        &\le \frac{1}{2} \norm{\nu \prt_2^2 \ozc^1
                - U^2 \prt_2 \ozc^1}_2  ^2
                + \frac{1}{2} \norm{w}_2^2, \\
    \abs{I_2}
        &\le \frac{\norm{\prt_2 U^2}_{L^\iny}}{2}
            \brac{\norm{\ozc^1}_2^2 + \norm{w}_2^2}.
\end{align*}
Both $I_1$ and $I_2$ vanish for constant $\UU$. For a varying $\UU$, $I_2$ can be easily controlled, since it only requires that $\norm{\ozc^1}_2$ vanish with $\nu$, which is required of any corrector. To control the first, the authors of \cite{GHT11} choose their corrector so it solves the elliptic equation (see (4.4) of \cite{GHT11}),
\begin{align}\label{e:GieHamoudaTemam2012Key}
	\nu \pdxn{\ozc^1}{x_2}{2} - U^2 \pdx{\ozc^1}{x_2} = 0,
\end{align}
which gives $I_2 = 0$. Applied to constant $U^2 = - U$, this leads to the expression for $\ozc^1$ in \cref{e:psiAndolz}.

The authors of \cite{GHT11} are led to \cref{e:GieHamoudaTemam2012Key} in a different manner, more in sympathy with \cite{TW02-1}. They suppose a scaling of the variables like Prandtl, though with a layer of width $\nu$ rather than $\sqrt{\nu}$, identifying in the resulting Prandtl-like expansion the key terms that control the expansion.

Finally, we note that the corrector as used by Kato in \cite{Kato1983} or see \cite{K2006Kato}---structureless in that it simply involves cutting off the stream function for the solution $\ou$ to the Euler equations---satisfies all of the bounds in \cref{e:zBounds}, and so requires one fewer derivative of regularity of $\ou$, and hence less involved compatibility conditions to obtain well-posedness of the Euler equations. Indeed, the well-posedness result of \cite{AKMRussianOriginal,AKM} would suffice to obtain the vanishing visocity limit using this simple corrector, except that the corrector's lack of structure makes it incapable of being used to bound the key term, $(\nu \prt_2 \ozc^1 - U^2 \ozc^1, \prt_2 w^1)$.

\appendix

\appendix
%
%
\section{Compatibility conditions}\label{A:CompConds}

\noindent We summarize in this appendix the compatibility conditions required to obtain a solution to $(E)$ as in \cref{P:GKM}. A more complete account is given in Section 4 of \cite{GKM2025}. The compatibility conditions are derived in \cite{GKM2025} for a 3D bounded domain, but are of the same form for a 2D bounded domain or periodic channel.

Recall that we use the notation $\uu^\BoldTau$ to be the tangential component of the velocity field $\uu$ along the boundary and we define the tangential vector field $\BoldTau$ so that $(\n, \BoldTau)$ is in the standard orientation of $(\e_1, \e_2)$. For our 2D channel, then, $\BoldTau = (-1, 0)$ and $\uu^\BoldTau = - u^1$ along the inflow boundary, while $\BoldTau = (1, 0)$ and $\uu^\BoldTau = u^1$ along the outflow boundary. 

Given $\uu$ with data regularity $N$ for some $N \ge 0$, we define the $N^{th}$ compatibility condition,
\begin{align}\label{e:condN}
	\begin{split}
	&\cond_{-1}:
		\ou_0^\BoldTau = \UU_0^\BoldTau
			\text{ on } \Gamma_+, \\
	&\cond_N: \cond_{N - 1} \text{ and } 
		\prt_t^{N + 1} U^\BoldTau|_{t = 0}
			= \wprt^{N + 1} \ou_0^\BoldTau
			\text{ on } \Gamma_+.
	\end{split}
\end{align}

For integers $n \ge 0$, we define $\wprt^n \ou_0$ inductively by setting $\wprt^0 \ou_0 = \ou_0$, while for $n \ge 1$, we take the time derivative of $\wprt^{n - 1} \uu$ at time zero and replace each instance of $\prt_t \uu$ in the resulting expression by $- \ou_0 \cdot \grad \ou_0 - \grad p^0 + \of(0)$. Here, $p^0$ is the value the pressure would have at time zero if $\uu$ actually solved $(E)$; that is, $p^0$ is the solution to
\begin{align}\label{e:p0}
	\begin{cases}
		\Delta p^0 = - \dv (\ou_0 \cdot \grad \ou_0)
			&\text{in } \Omega, \\
		\grad p^0 \cdot \n = -\prt_t U^\n(0) - \ou_0 \cdot \grad \ou_0
			&\text{on } \Gamma.
	\end{cases}
\end{align}

In other words, $\wprt^n \ou_0$ is the value that $\prt_t^n \uu$ would have at time zero if $\uu$ were a sufficiently regular solution to the Euler equations.

For $N = 0$, \cref{e:condN} is the compatibility condition in (1.10), (1.11) of Chapter 4 of \cite{AKM}:
\begin{align*}
	&\cond_0: 
		\prt_t U^\BoldTau|_{t = 0}
			= [- \ou_0 \cdot \grad \ou_0 - \grad p^0 + \of(0)]^\BoldTau
			\text{ on } \Gamma_+.	
\end{align*}

\bibliography{Refs}
\bibliographystyle{plain}

\end{document}